\baselineskip=15pt plus 2pt
\magnification=1150
\font\bigbold=cmbx7 scaled \magstep3
\noindent {\bigbold Transportation on spheres via an entropy formula}\par
\font\smallcaps=cmr8 scaled \magstep0
\vskip.05in
\indent {\bf Gordon Blower}\par
\indent {Department of Mathematics and Statistics,}\par
\indent {Lancaster University,}\par
\indent {Lancaster LA1 4YF ({\tt g.blower@lancaster.ac.uk})\par

\vskip.05in
\indent{13th July 2022}\par
\vskip.1in
{\narrower\smallskip\noindent The paper proves transportation inequalities for probability measures on spheres for the Wasserstein metrics with respect to cost functions that are powers of the geodesic distance. Let $\mu$ be a probability measure on the sphere ${\bf S}^n$ of the form $d\mu =e^{-U(x)}dx$ where $dx$ is the rotation invariant probability measure, and $(n-1)I+{\hbox{Hess}}\,U\geq {\kappa_U}I$, where $\kappa_U>0$. Then any probability measure $\nu$ of finite relative entropy with respect to $\mu$ satisfies ${\hbox{Ent}}(\nu\mid\mu) \geq (\kappa_U/2)W_2(\nu, \mu )^2$. The proof uses an explicit formula for the relative entropy which is also valid on connected and compact $C^\infty$ smooth Riemannian manifolds without boundary. A variation of this entropy formula gives the Lichn\'erowicz integral.\smallskip}
\vskip.1in
\indent {{\sl Key words:} Wasserstein metric; curvature; transport; convexity}\par
\vskip.1in
\indent 2020 {\sl Mathematics subject classification} 60E15, 58C35\par
\vskip.1in

\noindent {\bf 1. Transportation on the sphere}\par
\vskip.05in
\noindent Optimal transportation involves moving unit mass from one probability distribution to another, at minimal cost, where the cost is measured by Wasserstein's distance.\par
\vskip.05in

\noindent D{\smallcaps{EFINITION}} Let $(M,d)$ be a compact metric space and let $\mu$ and $\nu$ be probability measures on $M$. Then for $1\leq p<\infty$, Wasserstein's distance from $\mu$ to $\nu$ is $W_p(\nu, \mu )$, where
$$W_p(\nu, \mu )^p=\inf_\pi\Bigl\{ \int\!\!\!\int_{M\times M} d(x,y)^p \pi (dxdy): \pi\in {\hbox{Prob}}(M\times M)\Bigr\}\eqno(1.1)$$
where the probability measure $\pi$ has marginals $\nu$ and $\mu$. (See [8], [14].)\par
\vskip.05in
Transportation inequalities are results that bound the transportation cost $W_p(\nu, \mu )^p$ in terms of $\mu$, $\nu$ and geometrical quantities of $(M,d)$. Typically, one chooses $\mu$ to satisfy special conditions, and then one imposes minimal hypotheses on $\nu$. 
In this section, we consider the case where $(M,d)$ is the unit sphere ${\bf S}^2$ in ${\bf R}^3$, and obtain transportation inequalities by vector calculus. In section two, we extend these methods to a connected, compact and $C^\infty$ smooth Riemannian manifold $(M,d)$. \par

\indent On ${\bf S}^2$, let $\theta\in [0, 2\pi )$ be the longitude and $\phi\in [0, \pi]$ the colatitude, so the area measure is $dx=\sin\phi \, d\phi d\theta$. Let $ABC$ be a spherical triangle where $A$ is the North Pole; then by [10] the Green's function $G(B,C)=-(4\pi)^{-1}\log (1-\cos d(B,C))$ may be expressed in terms of longitude and co latitude of $B$ and $C$ via the spherical cosine formula. A related cost function is listed in [14], p 972.
Given probability measures $\mu$ and $\nu$ on ${\bf S}^2$, we can form 
$$G(\mu -\nu )(x)=\int_{{\bf S}^2}G(x,y)(\mu (dy)-\nu (dy))$$ 
with gradient in the $x$ variable 
$$\nabla G(\mu -\nu )(x)=\int_{{\bf S}^2}\nabla_xG(x,y)(\mu (dy)-\nu (dy)).$$
\par
\vskip.05in

\noindent {P{\smallcaps{ROPOSITION}} 1.1.} {\sl Let $\mu$ and $\nu$ be nonatomic probability measures on ${\bf S}^2$. Then 
$$ W_1(\mu,\nu )\leq \int_{{\bf S}^2}\Vert\nabla G(\mu-\nu )(x)\Vert dx.\eqno(1.2)$$} 
\vskip.05in
\noindent {\sl Proof.} The Green's function is chosen so that $\nabla\cdot \nabla G(B,C)=\delta_B(C)-1/(4\pi)$ in the sense of distributions. Given non-atomic probability measures $\mu$ and $\nu$ on ${\bf S}^2$, their difference $\mu-\nu$ is orthogonal to the constants on ${\bf S}^2,$ so for a $1$-Lipschitz function $\varphi: {\bf S}^2\rightarrow {\bf R}$, we have
$$\eqalignno{ \int_{{\bf S}^2} \varphi (x)(\mu (dx)-\nu (dx))&=\int_{{\bf S}^2} \varphi (x)\nabla\cdot\nabla G(\mu-\nu )(x) dx\cr
&=-\int_{{\bf S}^2} \nabla\varphi (x)\cdot\nabla G(\mu-\nu )(x)dx&(1.3)}$$
so by Kantorovich's duality theorem [8], the Wasserstein transportation distance is bounded by
$$ W_1(\mu,\nu )\leq \int_{{\bf S}^2}\Vert\nabla G(\mu-\nu )(x)\Vert dx.\eqno(1.4)$$ 
\vskip.1in

\noindent {D{\smallcaps{EFINITION}}} Suppose that $\mu$ is a probability measure and $\nu$ is a probability measure that is absolutely continuous with respect to $\mu$, so $d\nu =vd\mu$ for some probability density function $v\in L^1(\mu )$. Then the relative entropy of $\nu$ with respect to $\mu$ is 
 $${\hbox{Ent}}(\nu \mid\mu )=\int_{{\bf S}^2} \log v(y)\,\nu (dy),\eqno(1.5)$$
where $0\leq {\hbox{Ent}}(\nu \mid\mu )
\leq \infty$ by Jensen's inequality.\par
\vskip.05in
\indent At $x\in {\bf S}^2$, we have tangent space $T_s{\bf S}^2=\{ y\in {\bf R}^3: x\cdot y=0\}$. For $y\in T_x{\bf S}^2$ with $\Vert y\Vert=1$, we consider $\exp_x(ty)=x\cos t+y\sin t$ so that $\exp_x(0)=x$, $\Vert \exp_x(ty)\Vert=1$ and $(d/dt)_{t=0}\exp_x(ty)=y$; hence $\exp_x:T_x{\bf S}^2\rightarrow {\bf S}^2$ gives the exponential map. We let $J_{\exp_x}$ be the Jacobian determinant of this map. \par

\indent Suppose that $\mu (dx)=e^{-U (x)}dx$ is a probability measure and $\nu$ is a probability measure that is absolutely continuous with respect to $\mu$, so $d\nu =vd\mu$. We say that a Borel function $\Psi :{\bf S}^2\rightarrow {\bf S}^2$ induces $\nu$ from $\mu$ if $\int f(y)\nu (dy)=\int f(\Psi (x))\mu (dx )$ for all $f\in C({\bf S}^2; {\bf R})$. McCann [12] showed that there exists $\Psi$  that gives the optimal transport strategy for the $W_2$ metric; further, there exists a Lipschitz function $\psi: {\bf S}^2\rightarrow {\bf R}$ such that $\Psi (x)=\exp_x(\nabla\psi (x))$; so that
$$W_2(\nu, \mu )^2=\int_{{\bf S}^2} d(\Psi (x),x)^2 \mu (dx)= \int_{{\bf S}^2} \Vert \nabla\psi (x)\Vert^2 \mu (dx).\eqno(1.6)$$
\indent Talagrand developed $T_p$ inequalities in which $W_p(\nu, \mu )^p$ is bounded in terms of ${\hbox{Ent}}(\nu \mid\mu )$, as in [14] p 569.  In [5] and [6], the authors obtain some functional inequalities that are related to $T_p$ inequalities.  Here we offer an approach that is more direct, and uses only basic differential geometry to augment McCann's fundamental result. The key point is an explicit formula for the relative entropy in terms of the optimal transport maps.\par

\vskip.05in
\noindent {L{\smallcaps{EMMA}} 1.2.} {\sl Suppose that $\nu$ has finite relative entropy with respect to $\mu$, and let} 
$$H={\hbox{Hess}}_x\psi (x)\quad {{and}}\quad  A={\hbox{Hess}}_xd(x,y)^2/2\quad {{at}}\quad y=\Psi(x);\eqno(1.7)$$
{\sl let $\Psi_t(x)=\exp_x(t\nabla\psi(x))$ for $t\in [0,1]$. Then the relative entropy satisfies} 
$$\eqalignno{{\hbox{Ent}}(\nu \mid\mu )&\geq\int_{{\bf S}^2}\Bigl( {\hbox{trace}}\,\bigl( H-\log (A+H)\bigr)-\log J_{\exp_x}(\nabla\psi (x))\cr
&\quad +\int_0^1(1-t) {{d^2}\over{dt^2}}U (\Psi_t(x))dt\Bigr)\mu (dx).&(1.8)}$$
{\sl where $A$ is positive definite, $H$ is symmetric and $A+H$ is also positive definite, and}
$${\hbox{trace}}\,(H-\log (A+H))\geq 0.\eqno(1.9)$$
{\sl If $\psi\in C^2$, then equality holds in (1.8).} 
\vskip.1in
\noindent {\sl Proof.} To express the relative entropy in terms of the transportation map, we adapt an argument from [1]. We have 
${\hbox{Ent}}(\nu \mid\mu )=\int_{{\bf S}^2} \log v(\Psi (x))\mu (dx)$, where the integrand is
$$\log v(\Psi (x))=U (\Psi (x))-U (x)-\log J_\Psi (x),\eqno(1.10)$$
where the final term arises from the Jacobian of the change of variable $y=\Psi (x)$, where $\Psi=\Psi_1$ and $\Psi_t(x)=\exp_x(t\nabla\psi (x))$. We compute this Jacobian by the chain rule for derivatives with respect to $x$. Specifically by [6] p 622, we have ${\hbox{Hess}}(\psi (x)+d(x,y)^2/2)\geq 0$ and
$$\log J_\Psi (x)=\log J_{\exp_x}(\nabla \psi (x))+\log\det {\hbox{Hess}}(\psi (x)+d(x,y)^2/2)\eqno(1.11)$$
where $J_{\exp_x}$ is the Jacobian of $\exp_x:T_x{\bf S}^2\rightarrow {\bf S}^2$ and ${\hbox{Hess}}=D_x^2$ is the Hessian, where the expression is evaluated at $y=\exp_x(\nabla\psi (x))$. For $x\in {\bf S}^2$ and $\tau\in {\bf R}^3$ such that $x\cdot \tau =0$, we have $\tau\in T_x{\bf S}^2$ and 
$$\exp_x( \tau )=\cos (\Vert \tau\Vert )\, x+{{\sin (\Vert \tau\Vert )}\over{\Vert \tau \Vert}}\tau;\eqno(1.12)$$
see [5]. By a vector calculus computation, which we replicate from [5], one finds
$$J_{\exp_x}(\Vert\nabla \psi (x)\Vert) ={{\sin\Vert\nabla \psi (x)\Vert }\over{\Vert\nabla \psi (x)\Vert}}.\eqno(1.13)$$
 With $\psi :{\bf S}^2\rightarrow {\bf R}$ we have $\nabla\psi (x)\perp x$, so $0=x\cdot \nabla\psi (x),$ hence $0=\nabla\psi (x)+{\hbox{Hess}}(\psi (x)) x$. We write $\theta =\Vert\nabla \psi (x)\Vert$ for the angle between $x$ and $\Psi (x)$ so
$$\Psi (x)=\exp_x(\nabla \psi (x))=x\cos\theta +{{\sin\theta}\over{\theta}}\nabla\psi (x);$$
let $v=x\times \theta^{-1}\nabla\psi (x)$ where $\times$ denotes the usual vector product; then $\{ x, \theta^{-1}\nabla\psi (x), v\}$ gives an orthonormal basis of ${\bf R}^3$. Hence 
$${{\partial \Psi}\over{\partial  v}}=v\cos \theta -\sin\theta \langle \nabla\theta ,v\rangle x+\Bigl(\cos\theta -{{\sin\theta}\over{\theta}} \Bigr)\langle\nabla\theta ,v\rangle {{\nabla\psi (x)}\over{\theta}}+{{\sin\theta}\over{\theta}} {\hbox{Hess}}\psi (x) v,$$
\noindent and we obtain (1.13) from the final factor. Then by spherical trigonometry, we have
$$\cos d(\exp_x(\tau ), y)=(\cos \Vert\tau\Vert )\, \cos d(x,y)
+{{\sin \Vert\tau\Vert}\over{\Vert \tau\Vert}}\langle \tau ,y\rangle ,\eqno(1.14)$$
so we have $\langle \nabla_x \cos d(x,y), \tau \rangle=\langle y, \tau\rangle$ and $\langle {\hbox{Hess}}_x\cos d(x,y)\tau, \tau\rangle =-(\cos d(x,y)) \Vert\tau\Vert^2$;  so
$$\langle A\tau, \tau\rangle= {{1}\over{2}}\bigl\langle {\hbox{Hess}}_x d(x,y)^2\tau,\tau \bigr\rangle ={{d(x,y)}\over{\tan d(x,y)}} \Vert \tau\Vert^2 +\Bigl(1-{{d(x,y)}\over{\tan d(x,y)}}\Bigr) {{\langle y, \tau\rangle^2}\over{\sin^2d(x,y)}};\eqno(1.15)$$
hence $A$ is positive definite and is a rank-one perturbation of a multiple of the identity matrix. Note that the formulas degenerate on the cut locus $d(x,y)=\pi ;$ consider the international date line opposite the Greenwich meridian.\par
\indent We have
$${\hbox{Ent}}(\nu \mid\mu )=\int_{{\bf S}^2}\bigl( U (\Psi (x))-U (x)-\log J_\Psi (x)\bigr)e^{-U (x)}dx\eqno(1.16)$$
in which 
$$U (\Psi (x))-U (x)=\langle \nabla U (x), \nabla \psi (x)\rangle +\int_0^1 (1-t){{d^2}\over{dt^2}}U (\Psi_t(x))dt,\eqno(1.17)$$
and we can combine the first two terms in (1.16) by the divergence theorem
 so
$$\int_{{\bf S}^2} \langle \nabla U (x), \nabla\psi (x)\rangle e^{-U (x)} dx=\int_{{\bf S}^2} \nabla\cdot\nabla \psi (x) e^{-U (x)} dx.\eqno(1.18)$$
Hence from (1.11) we have
$${\hbox{Ent}}(\nu \mid\mu )= \int_{{\bf S}^2} \Bigl( \nabla\cdot \nabla \psi (x)-\log J_\Psi (x)\bigr) \mu (dx) +\int_0^1(1-t) {{d^2}\over{dt^2}}U (\Psi_t(x))dt\mu (dx),\eqno(1.19)$$
in which the Alexandrov Hessian [6],  [14]  p 363 satisfies
$${\hbox{trace}}\, {\hbox{Hess}}_x\psi (x)\leq\nabla\cdot\nabla \psi (x)=\Delta_D \psi (x),\eqno(1.20)$$
where $\Delta_D\psi$ is the distributional derivative of the Lipschitz function $\psi$; so we recognise (1.8).\par
\indent We have an orthonormal basis 
$$\Bigl\{ x, {{\nabla\psi (x)}\over{\Vert \nabla \psi (x)\Vert}}, x\times {{\nabla\psi (x)}\over{\Vert \nabla \psi (x)\Vert }}\Bigr\}\eqno(1.21)$$
for ${\bf R}^3$ in which the final two vectors give an orthonormal basis for $T_x{\bf S}^2$. Then 
$$\Bigl\langle A {{\nabla\psi (x)}\over{\Vert \nabla \psi (x)\Vert}},
{{\nabla\psi (x)}\over{\Vert \nabla \psi (x)\Vert }}\Bigr\rangle =1\eqno(1.22)$$
and 
$$\Bigl\langle A \Bigl(x\times {{\nabla\psi (x)}\over{\Vert \nabla \psi (x)\Vert}}\Bigr),
x\times {{\nabla\psi (x)}\over{\Vert \nabla \psi (x)\Vert}}\Bigr\rangle ={{d(x,y)}\over{\tan d(x,y)}},\eqno(1.23)$$
hence $A$ and $H$ have the form
$$A=\left[\matrix{1&0\cr 0& {{\Vert \nabla\psi (x)\Vert }\over{\tan \Vert\nabla\psi (x)\Vert}}\cr}\right], \qquad H=\left[\matrix{h&\beta \cr \beta &k\cr}\right]\eqno(1.24)$$
with respect to the stated basis of $T_x{\bf S}^2$.\par
\indent The function $f(x)=x-1-\log x$ for $x>0$ is convex and takes its minimum value at $f(1)=0$. Let $T$ be a self-adjoint matrix with eigenvalues $\lambda_1\geq \dots \geq\lambda_n$ where $\lambda_n>-1$; then the Carleman determinant of $I+T$ is $\det_2(I+T)=\prod_{j=1}^n (1+\lambda_j)e^{-\lambda_j}$. Since $A+H$ is positive definite, as in [1] Corollary 4.3, we can apply the spectral theorem to compute the Carleman determinant and show that 
$$-\log\det_2(A+H)={\hbox{trace}}\,\bigl( A+H-I-\log (A+H)\bigr)\geq 0\eqno(1.25)$$
so
$$\eqalignno{{\hbox{trace}}\,\bigl( H-\log (A+H)\bigr)&={\hbox{trace}}\,\bigl( A+H-I-\log (A+H)\bigr) +{\hbox{trace}}\, (I-A)\cr
&\geq 0+1-{{\Vert\nabla\psi (x)\Vert}\over{\tan\Vert\nabla\psi (x)\Vert}}\geq 0.&(1.26)}$$
\vskip.05in
\noindent {P{\smallcaps{ROPOSITION}} 1.3.} {\sl Suppose that the Hessian matrix of $U$ satisfies}
$${\hbox{Hess}}\,U (x)+ I\geq \kappa_U I\qquad (x\in {\bf S}^2)\eqno(1.27)$$
{\sl for some $\kappa_U>0$. Then $\mu$ satisfies the transportation inequality}
$${\hbox{Ent}}(\nu\mid\mu )\geq {{\kappa_U}\over{2}}W_2(\nu ,\mu )^2.\eqno(1.28)$$
{\sl This applies in particular when $\mu$ is normalized surface area measure.}\par
\vskip.1in
\noindent {\sl Proof.} Let $K:[0, \pi )\rightarrow {\bf R}$ be the function
$$K(\alpha )=1-{{\alpha}\over{\tan\alpha}}+\log {{\alpha}\over{\sin\alpha}}={{d}\over{d\alpha}}\Bigl( \alpha \log {{\alpha}\over{\sin\alpha}}\Bigr).\eqno(1.29)$$
Then from (1.13) and (1.26) we have

$$\int_{{\bf S}^2} \Bigl( \nabla\cdot \nabla \psi (x)-\log J_\Psi (x)\bigr) \mu (dx)\geq \int_{{\bf S}^2} \Bigl( -\log\det_2(A+H)+K(\Vert\nabla \psi (x)\Vert)\Bigr)\mu (dx).$$ 

 Considering the final integral in (1.8), we have 
$${{\partial\Psi_t(x) }\over{\partial t}}=-\Vert\nabla\psi(x)\Vert\sin (t\Vert\nabla\psi (x)\Vert) x+\cos (t\Vert\nabla\psi (x)\Vert)\nabla\psi (x)\eqno(1.30)$$
which has constant speed $\Vert {{\partial\Psi_t(x) }\over{\partial t}}\Vert=\Vert\nabla \psi (x)\Vert$ and $\langle {{\partial\Psi_t(x) }\over{\partial t}}, \Psi_t(x)\rangle=0;$ also
$${{\partial^2}\over{\partial t^2}}U (\Psi_t(x))=\Bigl\langle {\hbox{Hess}}U \circ \Psi_t(x) {{\partial\Psi_t(x) }\over{\partial t}}, {{\partial\Psi_t(x) }\over{\partial t}}\Bigr\rangle -\Vert\nabla\psi (x)\Vert^2 \bigl\langle (\nabla U )\circ\Psi_t(x), \Psi_t(x)\bigr\rangle,\eqno(1.31)$$
where the final term is zero since $\nabla U\circ\Psi_t(x)$ is in the tangent space at $\Psi_t(x)$, hence is perpendicular to $\Psi_t(x)$.
We therefore have the crucial inequality
$$\eqalignno{{\hbox{Ent}}(\nu\mid\mu )&\geq\int_{{\bf S}^2} \Bigl( -\log\det_2(A+H)+K(\Vert\nabla \psi (x)\Vert)\cr
&\quad +\int_0^1(1-t)\Bigl\langle {\hbox{Hess}}U\circ \Psi_t(x){{\partial\Psi_t(x) }\over{\partial t}}, {{\partial\Psi_t(x) }\over{\partial t}}\Bigr\rangle dt\Bigr)\mu (dx)&(1.32)}$$
To simplify the function $K$,  we recall from [9] 8.342 the Maclaurin series
$$\eqalignno{\log {{\alpha}\over{\sin\alpha}}&=\log\Gamma \Bigl(1+ {{\alpha}\over{\pi}}\Bigr)+\log\Gamma \Bigl(1-{{\alpha}\over{\pi}}\Bigr)\cr
&=\sum_{m=1}^\infty {{\zeta (2m)}\over{\pi^{2m}m}}\alpha^{2m}\qquad (\vert\alpha\vert<\pi ),&(1.33)}$$
where we have introduced Euler's $\Gamma$ function and Riemann's $\zeta$ function, so
$$\eqalignno{K(\alpha )&=\sum_{m=1}^\infty {{(2m+1)\zeta (2m)}\over{\pi^{2m}m}}\alpha^{2m}\geq {{3\zeta(2)}\over{\pi^2}}\alpha^2={{\alpha^2}\over{2}}.&(1.34)}$$
Now we consider (1.32) with the hypothesis (1.27) in force. The Carleman determinant contributes a nonnegative term as in (1.25), while the final integral in (1.32) combines with the integral of $K(\Vert \nabla\psi (x)\Vert)$ to give
$$\eqalignno{{\hbox{Ent}}(\nu \mid \mu )&\geq \int_{{\bf S}^2} \bigl( K(\Vert \nabla\psi (x)\Vert )+(1/2)\Vert \nabla\psi (x)\Vert^2\bigr)\mu (dx)\cr
&\geq{{\kappa_U}\over{2}}\int_{{\bf S}^2}\Vert \nabla\psi (x)\Vert^2\mu (dx)\cr
&={{\kappa_U}\over{2}}W_2(\nu,\mu)^2.&(1.35)}$$ 
\indent When $\mu$ is normalized surface area, $U$ is a constant and the hypothesis (1.27) holds with $\kappa_U=1$.
\vskip.05in
\vskip.05in
\noindent {\bf 2. Transportation on compact Riemannian manifolds}\par
\vskip.05in
\noindent Let $M$ be a connected, compact and $C^\infty$ smooth Riemannian manifold of dimension $n$ without boundary, and let $g$ be the Riemannian metric tensor, giving metric $d$. Let $\mu (dx)=e^{-U(x)}dx$ be  a probability measure on $M$ where $dx$ is Riemannian measure and $U\in C^2(M; {\bf R})$. 
Suppose that $\nu$ is a probability measure on $M$ that is of finite relative entropy with respect to $\mu$. Then by McCann's theory [12], there exists a Lipschitz function $\psi :M\rightarrow {\bf R}$ such that $\Psi (x)=\exp_x(\nabla\psi (x))$ induces $\nu$ from $\mu$.  then we let $\Psi_t(x)=\exp_x(t\nabla\psi (x))$. We proceed to compute quantities which we need for our extension of Lemma 1.2.\par
\indent 
Given distinct points $x,y\in M$, we suppose that $x=\exp_y(\xi )$, and for $w\in T_yM$ introduce
$$\gamma (s,t)= \exp_y(t(\xi +sw))\eqno(2.1)$$
so that $t\mapsto \gamma (s,t)$ is a geodesic, and in particular $\gamma (0,t)$ is the geodesic from $y=\gamma (0,0)$ to $x=\gamma (0,1)$.  When $y=\exp_x(\nabla\psi (x))$ for a Lipschitz function $\psi :M\rightarrow {\bf R}$, we can determine $\xi$ as follows. Let $\phi (z)=-\psi (z)$ and introduce its infimal convolution
$$\phi^c (y)=\inf_w\{ (1/2)d(y,w)^2-\phi (w)\}\eqno(2.2)$$
which is attained at $x$ since $y=\exp_x(\nabla\psi (x))=\exp_x(-\nabla \phi (x))$. Now $\phi^{cc}(x)=\phi (x)$, so
$$\phi (x)=\inf_w\{ (1/2)d(x,w)^2-\phi^c (w)\}\eqno(2.3)$$
where the infimum is attained at $y$ since $\phi (x)+\phi^c(y)=d(x,y)^2/2$. By Lemma 2 of [12], $\phi^c$ is Lipschitz and
$$x=\exp_y(-\nabla \phi^c(y)).\eqno(2.4)$$
The speed of $\gamma (0,t)$ is given by
$$\Bigl\Vert {{\partial\gamma}\over{\partial t}}\Bigr\Vert=\Vert\nabla \phi^c(y)\Vert =d(y, \exp_y(-\nabla\phi^c(y)))=d(x,y)=d(x,\exp_x(-\nabla\phi (x)))=\Vert\nabla\psi (x)\Vert .\eqno(2.5)$$
\indent Let $R$ be the curvature of the Levi-Civita derivation $\nabla$ so $$R(X,Y)Z=\nabla_Y\nabla_XZ-\nabla_X\nabla_YZ-\nabla_{[Y,X]}Z\qquad (X,Y,Z\in T_xM).$$
Then by [13] page 36,  for all $Y\in T_xM$, the curvature operator $R_Y: X\mapsto R(X,Y)Y$ is self-adjoint with respect to the scalar product on $T_xM$. Also 
$$Y(s,t)={{\partial}\over{\partial s}}\gamma (s,t)\eqno(2.6)$$
satisfies the initial conditions
$$Y(s,0)=0,\quad {{\partial Y}\over{\partial t}}(0,0)=w,\eqno(2.7)$$
and Jacobi's differential equation [4] (2.43)
$${{\partial^2Y}\over{\partial t^2}}+R\Bigl( {{\partial\gamma}\over{\partial t}}, Y\Bigr){{\partial\gamma}\over{\partial t}} =0.\eqno(2.8)$$
 By calculating the first variation of the length formula [13] p 161, one shows that
$${{1}\over{2}}\Bigl\langle {\hbox{Hess}}_x d(x,y)^2 Y(0,1),Y(0,1)\Bigr\rangle =g\Bigl({{\partial Y}\over{\partial t}}(0,1), Y(0,1)\Bigr).\eqno(2.9)$$
Assume that there are no conjugate points on $\gamma (s,t)$. Then by varying $w$, we can make $Y(0,1)$ cover a neighbourhood of $0$ in $T_xM$.
 Let 
$$A={{1}\over{2}}{\hbox{Hess}}_x d(x,y)^2\Bigr\vert_{y=\exp_x(\nabla\psi (x))},\eqno(2.10)$$
and 
$$H ={\hbox{Hess}} \psi (x).\eqno(2.11)$$
Let $J_{\exp_x}(v)$ be the Jacobian of the map $T_xM\rightarrow M$ given by $v\mapsto \exp_x(v)$, as in (3.4) of [3].  
\vskip.05in
\noindent {L{\smallcaps{EMMA}} 2.1.} {\sl  Suppose that $\Psi_t (x)=\exp_x(t \nabla\psi (x))$, where $\Psi_1$ induces the probability measure $\nu$ from $\mu$ and gives the optimal transport map for the $W_2$ metric. Then the relative entropy satisfies} 
$$\eqalignno{{\hbox{Ent}}(\nu \mid\mu )&\geq\int_{M}\Bigl( {\hbox{trace}}\,\bigl( H-\log (A+H)\bigr)-\log J_{\exp_x}(\nabla\psi (x))\cr
&\quad +\int_0^1(1-t) \Bigl\langle {\hbox{Hess}}\,U\circ \Psi_t(x){{\partial \Psi_t(x)}\over{\partial t}}, {{\partial\Psi_t(x)}\over{\partial t}}\Bigr\rangle dt\Bigr)\mu (dx).&(2.12)}$$
{\sl where $H$ is symmetric and $A+H$ is also positive definite. If $\psi\in C^2(M; {\bf R})$, then equality holds in (2.12).}\par
\vskip.1in
\noindent {\sl Proof.} This is similar to Lemma 1.2. As in (125),  we have
$$\eqalignno{{\hbox{trace}}\,\bigl( H-\log (A+H)\bigr)&=-\log\det_2(A+H)+{\hbox{trace}}(I-A)\cr
&\geq {\hbox{trace}}(I-A),&(2.13)}$$
and by standard calculations [13] p32 we have
$${{\partial^2}\over{\partial t^2}}U (\Psi_t(x))=\Bigl\langle {\hbox{Hess}}U\circ \Psi_t(x){{\partial \Psi_t(x)}\over{\partial t}}, {{\partial\Psi_t(x)}\over{\partial t}}\Bigr\rangle\eqno(2.14)$$
since $\Psi_t(x)$ is a geodesic. \par
\vskip.05in
\indent The curvature operator is the symmetic operator $R_Z:Y\mapsto R(Z,Y)Z$. If $M$ has nonnegative Ricci curvature so that $R_Z\geq 0$ as a matrix for all $Z$, then we have
$$-\log J_{\exp_x}(\nabla\psi (x))\geq 0.\eqno(2.15)$$
by (3.4) of [Ca].\par

The following result recovers the Lichn\'erowicz integral, as in (4.16) of [1] and (1.1) of [7]. This integral also appears implicitly in the Hessian calculations in Appendix D of [11]. Let $\Vert H\Vert_{HS}$ be the Hilbert--Schmidt norm of $H$.\par
\vskip.05in

\noindent {P{\smallcaps{ROPOSITION}} 2.2.} {\sl  Suppose that $\psi\in C^2(M; {\bf R})$ and $\Psi_\tau (x)=\exp_x(\tau \nabla\psi (x))$ induces a probability measure $\nu_\tau$ from $\mu$ such that $\Psi_\tau$ is the optimal transport map for the $W_2$ metric. Then}
$$\eqalignno{{\hbox{Ent}}(\nu_\tau \mid \mu )&={{\tau^2}\over{2}}\int_{M} \Bigl( \Vert {\hbox {Hess}}\,\psi (x)\Vert_{HS}^2 +{\hbox{trace}}\, R_{\nabla\psi (x)}\cr
&\qquad +\bigl\langle {\hbox{Hess}}\,U (x)\nabla\psi (x),\nabla\psi (x)\bigr\rangle\Bigr) \mu (dx)+O(\tau^3)\qquad (\tau\rightarrow 0+).&(2.16)\cr}$$
\vskip.05in
\noindent {\sl Proof.} For small $\tau>0$, we rescale $\psi$ to $\tau\psi$ and consider $y=\exp_x(\tau\nabla\psi (x))$; then we return to $x$ along a geodesic $\gamma_\tau (t)=\exp_y(-t\nabla (-\tau\psi )^c(y))$ for $0\leq t\leq 1$ with constant speed $\tau\Vert\nabla\psi (x)\Vert$. 
Observe that $\tau\psi (x)=(-\tau\psi)^c(y)-\tau^2\Vert\nabla\psi (x)\Vert^2/2$, and $\nabla_xd(x,y)^2/2=-\exp_x^{-1}(y)=-\tau\nabla \psi (x)$ and $\nabla_yd(x,y)^2/2=-\exp_y^{-1}(x)=\nabla (-\tau\psi )^c(y)$ by Gauss's Lemma.  Recalling that the curvature operator is self-adjoint by page 36 of [13], we choose the basis of $T_yM$ so that the first basis vector points along the direction of the geodesic 
$\gamma_\tau (0)$. Hence Jacobi's equation (2.8) can be expressed as a second order differential equation in block matrix form, with a symmetric matrix $ S_{-\nabla(-\tau\psi )^c(y)} $ given by components of the curvature tensor such that
$$R\Bigl( {{d\gamma_\tau}\over{dt}},Y\Bigr){{d\gamma_\tau}\over{dt}} =\left[\matrix{0&0\cr 0&S_{-\nabla(-\tau\psi )^c(y)}}\right]Y\qquad (0<t<1).\eqno(2.17)$$
as in (2.4) of [6]. Then the Jacobi equation reduces to a first-order block matrix equation with blocks of shape $(1+(n-1))\times (1+(n-1))$ in a $(2n)\times (2n)$ matrix
$${{d}\over{dt}}\left[\matrix{Y\cr V}\right] =\left[\matrix{0&0& 1&0\cr 0&0&0&I_{n-1}\cr 0&0&0&0\cr 0&-S_{-\nabla (-\tau\psi )^c(y)}&0&0}\right] \left[\matrix{Y\cr V}\right] ;\quad 
\left[\matrix{Y(0)\cr V(0)}\right] =\left[\matrix{0\cr w}\right].\eqno(2.18)$$

\vskip.05in

 To find the limit as  $\tau\rightarrow 0$, we can assume that $S_{-\nabla (-\tau\psi )^c (y)}$ is constant on the geodesic, and may be expressed as $\tau^2 S$  where $\tau^2 S=S_{\tau\nabla\psi (x)}$ has shape $(n-1)\times (n-1)$.
The functions $\cos \alpha$ and $\sin \alpha /\alpha$ are entire and even, so $\cos \sqrt {s}$ and $\sin \sqrt{s}/\sqrt{s}$ are entire functions, hence they operate on complex matrices. Note that the matrix
$$T=\left[\matrix{0&0\cr 0&S_{-\nabla(-\tau\psi)^c(y)}\cr}\right]$$
 in the bottom left corner is symmetric, has rank less than or equal to $n-1$, and does not depend upon $t$. Hence we consider the matrix  
$$\left[\matrix{Y\cr V}\right] =\left[\matrix{\cos (t\sqrt {T})& {{\sin( t\sqrt {T})}\over{\sqrt {T}}}\cr -\sqrt T\sin (t\sqrt{ T})& \cos ( t\sqrt { T})}\right]
\left[\matrix{Y_0\cr V_0}\right]$$
which has derivative
$${{d}\over{dt}}\left[\matrix{Y\cr V}\right] =\left[ \matrix{0&I\cr -T&0\cr}\right]\left[\matrix{\cos (t\sqrt {T})& {{\sin (t\sqrt {T})}\over{\sqrt {T}}}\cr -\sqrt {T}\sin (t\sqrt {T})& \cos (t\sqrt {T})}\right]\left[\matrix{Y_0\cr V_0}\right]$$
so we can use this formula to solve (2.18).
 So the approximate differential equation has solution 
$$\left[\matrix{Y(1)\cr V(1)}\right] =\left[ \matrix{ 1&0&1&0\cr 0&\cos \tau\sqrt{S} &0& {{\sin \tau\sqrt{S}}\over{\tau\sqrt{S}}}\cr
0&0&1&0\cr 0&-\tau\sqrt{S}\sin \tau\sqrt{S}&0&\cos \tau\, \sqrt{S}\cr}\right]\left[\matrix{0\cr w}\right].\eqno(2.19)$$
 Hence by (2.9) we have
$$A=\left[\matrix{ 1&0\cr 0&{{\tau\sqrt{S}}\over{\tan\tau\sqrt{S}}}\cr}\right]=(1+O(\tau^2))I_n\eqno(2.20)$$
which gives rise to the approximation
$${\hbox{trace}}(I_n-A)={\hbox{trace}}\Bigl( I_{n-1}-{{\tau\sqrt{S}}\over{\tan\tau\sqrt{S}}}\Bigr)= 
{{\tau^2}\over{3}}{\hbox{trace}}(S)+O(\tau^4)\qquad (\tau\rightarrow 0+),\eqno(2.21)$$ 
and likewise we obtain
$$-\log J_{\exp_x}(\tau \nabla\psi (x)) =-\log\det {{\sin \tau\sqrt{S}}\over{\tau\sqrt{S}}}={{\tau^2}\over{6}}{\hbox{trace}}(S)+O(\tau^4).\eqno(2.22)$$
From (2.19), we have
$$\eqalignno{-\log\det_2 (A+\tau H)&={{1}\over{2}}{\hbox{trace}}\bigl( (A-I_n+\tau H)^2\bigr)+O(\tau^3)\cr
&={{\tau^2}\over{2}}{\hbox{trace}}(H^2)+O(\tau^3)\cr
&={{\tau^2}\over{2}}\Vert {\hbox{Hess}}\, \psi (x)\Vert_{HS}^2+O(\tau^3),&(2.23)\cr}$$
so the result follows by Lemma 2.1. \par
\vskip.05in
We conclude with a transportation inequality which generalizes Proposition 1.3 to the unit spheres ${\bf S}^n$. See [2] for a discussion of measures on product spaces.\par
\vskip.05in
\noindent {T{\smallcaps{HEOREM}} 2.3.} {\sl Let $M={\bf S}^n$ for some $n\geq 2$,  and suppose that} 
$$(n-1)I+{\hbox{Hess}} \,U (x)\geq \kappa_U I\qquad (x\in {\bf S}^{n})\eqno(2.24)$$
{\sl for some $\kappa_U>0$. Then} 
$${\hbox{Ent}}(\nu\mid\mu )\geq {{\kappa_U}\over {2}}W_2(\nu, \mu )^2.\eqno(2.25)$$
\vskip.05in
\noindent {\sl Proof.} In this case, the curvature operator is constant, so we have $S_{\nabla\psi (x)}Y=\Vert \nabla \psi (x)\Vert^2Y$, so 
$${\hbox{trace}}\,R_{\nabla\psi (x)} =(n-1)\Vert \nabla\psi (x)\Vert^2.\eqno(2.26)$$
Thus the result follows with a similar proof to Proposition 1.3 using data from the proof of Proposition 2.2.\par

\vskip.05in
\noindent {\bf Acknowledgement} I thank Graham Jameson for helpful remarks concerning inequalities which led to (1.34). I am also grateful to the referee, whose helpful comments improved the exposition.\par
\vskip.05in
\noindent  {\bf References}\par
\vskip.05in
\noindent [1] G. Blower, The Gaussian isoperimetric inequality and transportation, {\sl Positivity} {\bf 7} (2003), 203-224.\par
\noindent [2] G. Blower and F. Bolley, Concentration of measure on product spaces with applications to Markov processes, {\sl Studia Math.} {\bf 175} (2006), 47-72.\par
\noindent [3] X. Cabre, Nondivergent elliptic equations on manifolds with nonnegative curvature, {\sl Comm. Pure Appl. Math.} {\bf  50} (1997), 623-665.\par
\noindent [4] I. Chavel, {\sl Riemannian Geometry: a modern introduction}, (Cambridge University Press, 1993).\par 
\noindent [5] D. Cordero-Erausquin, Pr\'ekopa--Leindler inequalities sur la sph\`ere, {\sl C.R. Acad. Sci. Paris} {\bf 329} (1999), 789-792.\par
\noindent [6] D. Cordero-Erausquin, R.J. McCann and M. Schmuckensl\"ager, Pr\'ekopa--Leindler type inequalities on Riemannian manifolds, Jacobi fields and optimal transport, {\sl Annales. de la Fac. Sci. Toulouse Math.} {\bf 15} (2006), 613-635.\par
\noindent [7] J.-D. Deuschel and D.W. Stroock, Hypercontractivity and spectral gap of symmetric diffusions with applications to stochastic Ising models, {\sl J. Funct. Anal.} {\bf 92} (1990), 30-48.\par
\noindent [8] R. M. Dudley, {\sl Real Analysis and Probability,} second edition, (Cambridge University Press, 2004).\par
\noindent [9] I. S. Gradsteyn and I.M. Ryzhik, {\sl Table of Integrals, Series and Products}, (Academic Press, 1965).\par
\noindent [10] Y. Kimura and H. Okamoto, Vortex motion on a sphere, {\sl J. Phys. Soc. Japan} {\bf 56} (1987), 4203-4206.\par
\noindent [11] J. Lott and C. Villani, Ricci curvature for metric measure spaces via optimal transport, {\sl Annals of Math.} (2) {\bf 169} (2009), 903-991.\par 
\noindent [12] R. J. McCann, Polar factorization of maps on Riemannian manifolds, {\sl Geom. Funct. Anal.} {\bf 11} (2001), 589-608.\par 
\noindent [13] P. Pedersen, {\sl Riemannian Geometry}, second edition, (Springer, 2006).\par
\noindent [14] C. Villani, {\sl Optimal Transport: Old and New}, (Springer, 2009).\par
\vfill
\eject
\end

For $y\in{\bf S}^n$,  we choose the basis of $T_y{\bf S}^2$ so that the first basis vector points along the direction of the geodesic $\gamma (0,t)$, which has constant speed $\rho=\Vert\nabla\psi (x)\Vert>0$ by (2.5), so that 
$$R=\left[\matrix{0&0\cr 0&\rho^2 I_{n-1}}\right].\eqno(2.16)$$
Then the Jacobi equation reduces to a first-order block matrix equation with blocks $(1+(n-1))\times (1+(n-1))$ in a $(2n)\times (2n)$ matrix
$${{d}\over{dt}}\left[\matrix{Y\cr V}\right] =\left[\matrix{0&0& 1&0\cr 0&0&0&I_{n-1}\cr 0&0&0&0\cr 0&-\rho^2I_{n-1}&0&0}\right] {{d}\over{dt}}\left[\matrix{Y\cr V}\right] ; 
\left[\matrix{Y(0)\cr V(0)}\right] =\left[\matrix{0\cr w}\right],\eqno(2.17)$$

with solution 
$$\left[\matrix{Y\cr V}\right] =\left[ \matrix{ 1&0&t&0\cr 0&\cos \rho t\, I_{n-1} &0& {{\sin \rho t}\over{\rho}}I_{n-1}\cr
0&0&1&0\cr 0&\rho\sin \rho t\, I_{n-1}&0&\cos \rho t\, I_{n-1}\cr}\right]\left[\matrix{0\cr w}\right]\eqno(2.18)$$
where $\rho =\Vert \nabla \psi (x)\Vert$. Hence
$$A=\left[\matrix{ 1&0\cr 0&{{\rho}\over{\tan\rho}}I_{n-1}}\right]=I_n+O(\tau^2)\eqno(2.19)$$
and 
$${\hbox{trace}}(I_n-A)=(n-1)\Bigl( 1-{{\rho}\over{\tan\rho}}\Bigr)\geq (n-1) {{\rho^2}\over{3}}\qquad (0<\rho <\pi),\eqno(2.20)$$ 
where the final step follows from elementary calculus applied to 
$$f(\rho )=-\rho\cos\rho+\sin\rho-(1/3)\rho^2\sin\rho.\eqno(2.21)$$

By the Lemma, we have
$${\hbox{Ent}}(\nu\mid\mu )\geq \int_{{\bf S}^n} \Bigl( {{(n-1)\Vert\nabla\psi (x)\Vert^2}\over{3}}+\int_0^1(1-t)
\Bigl\langle{\hbox{Hess}}U\circ \Psi_t(x){{\partial\Psi_t(x)}\over{\partial t}}, {{\partial\Psi_t(x)}\over{\partial t}}\Bigr\rangle dt\Bigr) \mu (dx)\eqno(2.22)$$
where $\Vert {{\partial\Psi_t(x)}\over{\partial t}}\Vert=\Vert\nabla\Psi (x)\Vert,$ so 
$${\hbox{Ent}}(\nu\mid\mu )\geq {{\kappa_U}{2}}\int_{{\bf S}^n}\Vert \nabla\psi (x)\Vert^2\mu (dx)={{\kappa_U}\over{2}}W_2(\nu, \mu )^2.\eqno(2.23)$$

Given a $1$-Lipschitz function $\varphi :{\bf S}^2\rightarrow {\bf R}$ and a Lipschitz function $\psi :{\bf S}^2\rightarrow {\bf R}$, we can introduce the one-parameter family of functions $\Psi_t: {\bf S}^2\rightarrow {\bf S}^2$ by $\Psi_t(x)=\exp_x (t\nabla \psi (x))$. Then 
$$\eqalignno{t\Vert \nabla \psi (x)\Vert &=d(\exp_x (t\nabla \psi (x)), x)\cr
&\geq \varphi (\Psi_t(x))-\varphi (x)\cr
&=t\langle \nabla \varphi (x), \nabla\psi (x)\rangle +(t^2/2)\langle {\hbox{Hess}} \, \varphi (x), \nabla \psi (x), \nabla \psi (x)\rangle +o(t^2)&(1.7)}$$
as $t\rightarrow 0+$, where ${\hbox{Hess}}\, \varphi $ is the Hessian in the sense of Alexandrov. We can take $\nabla \varphi (x)=\nabla \psi (x)/\Vert\nabla\psi (x)\Vert$ to maximize the right-hand side. \par